# A Robust Scientific Machine Learning for Optimization:

# A Novel Robustness Theorem

Luana P. Queiroz[1,2], Carine M. Rebello[3,*], Erber A. Costa[3], Vinicius V. Santana[1,2], Alírio E. Rodrigues[1,2], Ana M. Ribeiro[1,2] and Idelfonso B. R. Nogueira[1,2,*]

[1] Laboratory of Separation and Reaction Engineering, Associate Laboratory LSRE/LCM, Department of Chemical Engineering, Faculty of Engineering, University of Porto, Rua Dr. Roberto Frias, 4200-465 Porto, Portugal.

[2] ALiCE—Associate Laboratory in Chemical Engineering, Faculty of Engineering, University of Porto, Rua Dr. Roberto Frias, 4200-465 Porto, Portugal.

[3] Departamento de Engenharia Química, Escola Politécnica (Polytechnic Institute), Universidade Federal da Bahia, Salvador 40210-630, Brazil.

* Correspondence: carine.menezes@ufba.br; idelfonso@fe.up.pt

# Abstract

Scientific machine learning (SciML) is a field of increasing interest in several different application fields. In an optimization context, SciML-based tools have enabled the development of more efficient optimization methods. However, implementing SciML tools for optimization must be rigorously evaluated and performed with caution. This work proposes the deductions of a robustness test that guarantees the robustness of multiobjective SciML-based optimization by showing that its results respect the universal approximator theorem. The test is applied in the framework of a novel methodology which is evaluated in a series of benchmarks illustrating its consistency. Moreover, the proposed methodology results are compared with feasible regions of rigorous optimization, which requires a significantly higher computational effort. Hence, this work provides a robustness test for guaranteed robustness in applying SciML tools in multiobjective optimization with lower computational effort than the existent alternative.

## 1. Introduction

Classically in several engineering fields, phenomenological models are employed to represent the actual state of a given system. These models are accurate and allow extrapolations. However, real-time information is demanded in the current industry and processes scenario. Although the phenomenological models are robust, they require heavy computational resources, which leads to not being able to fulfill the requirement of providing online information. This need is found in many applications such as control, optimization, and inference. In this context, in the past decades, with prominence in the last decade, simulation-based and data-driven systems have been developed and rising in popularity due to the achievement of more powerful and faster computers [1]. These advanced models have been labor-saving techniques to overcome the challenges of optimization systems with the development of artificial intelligence-based surrogate models.

Still, the surrogate model is an approximation of the actual system. These models do not have a predictive capacity as the phenomenological models. Furthermore, they have limited ability to perform extrapolation. Another downside is that they can present artificial minimums. So, techniques that consider these problems must be developed to apply this strategy properly. However, this care is not usually observed in the literature.

For example, Li et al. (2017) [2] applied surrogate-based models to optimize orifice plates to improve the performance in water treatment. The systems' performance was analyzed by degradation quality and processing capacity per unit time. Although the results obtained solved the authors' problem efficiently, the optimality of the SciML models was not evaluated. In other words, the optimum point acquired was not proven to represent the system properly. Ochoa-Estopier et al. (2018) [3] optimized the product yield and debottlenecking the distillation units of a refinery with a surrogate-based model. This optimization achieved a profit increase of $7.2 million per year by manipulating the crude oil flow rate, furnace outlet temperatures, and stripping steam. However, the optimality of the result was not evaluated [3], leaving room to speculate that the profit increase could be higher than actually observed. Ho et al. (2019) [4] applied the same method as the last ones to determine the optimal productivity for biocatalytic continuous manufacturing of sitagliptin, an antidiabetic drug. The productivity optimization resulted in a significantly better performance of the system. Nonetheless, the optimality of this result also was not studied.

On the other hand, Nogueira et al. (2022) [5] proposed a methodology for concomitantly material screening and process optimization by applying surrogate-based optimization. The work concluded that the system studied can have a series of materials that can be implemented to obtain an optimal performance subjected to the operating conditions of the process. The authors verified if the optimal solution obtained was following the phenomenological model. Tests were done to uphold the nonexistence of artificial minima that could have been introduced by the surrogate-based model. However, these tests were limited to evaluating a few optimal points obtained. Pai et al. (2020) [6] optimized the process conditions for an adsorption-based process, also applying surrogate-based optimization for a $CO_2$ capture system. Their goal was to optimize the process energy consumption and productivity for a minimum process purity of 95 % and recovery of 80 %. The authors also evaluated the optimal point obtained compared to the phenomenological model. The optimal value provided by the surro-gate-based optimization differed by about 3 % from the value given by the detailed model. The optimal result was considered valid. A few more reports of surrogate-based models of optimization that had the optimal result evaluated and deemed accurate by comparison with the phenomenological model result were found in the literature, such as the work of Subraveti et al. (2019) [7] and Pai et al. (2020) [8]. Still, the authors mentioned above just evaluated the optimal point obtained and not the optimality of the overall optimization problem. Neither was assessed the consistency of the surrogate model identification to face an optimization context.

In this context, this work proposes the deductions of a test that assures that the optimization result based on SciML respects the universal approximation theorem. Hence, assessing the robustness of the optimal result obtained from a surrogate-based optimization. Therefore, certifying that the overall optimization result is reliable and can be confidently employed. This methodology is based on a meta-heuristic method that uses a population to perform the optimization. Then, using the population history and a statistical sampling method, it is proposed to evaluate the robustness of the results using the deduced test. This is done randomly by sampling the optimization population without replacement and assuring that there is no correlation between samples by the design of experiments. The samples are then evaluated in a robust source of information, which can be either a rigorous model or a real system. Then the deduced robustness metric is computed to evaluate how precisely the surrogate model outputs are following the rigorous system. In this way, the proposed methodology leverages the potential of a surrogate-based optimization while providing a robust assessment of its results, using few system evaluations in a more rigorous source of information.

The methodology proposed is evaluated using benchmark functions, providing a practical demonstration of the proposed strategy. Furthermore, the results are compared to the scenarios where no robustness assessment is done, which is the usual practice in the literature. The benefits of this evaluation are highlighted.

## 2. Materials and Methods

The overall description of the proposed methodology is depicted in Figure 1. It starts with a system of interest, which can be represented in different ways, such as through surrogate, phenomenological models, or historical data. Then, the training data acquisition is performed. This is an essential step that carries a heavyweight throughout all the optimization as the nature of the acquired data will define the model's success. The surrogate models are data-oriented, extracting information and learning from them. In case of inequality or problems in this phase, the model will learn those problems found in the dataset. Therefore, the methodology for data acquisition must be consistent and reliable. Moreover, the amount of data provided must be adequate so the model can identify patterns, learn from them, and be validated and tested.

Once the data is acquired, it is necessary to identify the surrogate model, which is also a crucial step that requires consistent methods and criteria to evaluate the performance of the models and define the best structure to be applied along with the next steps. After the model is well structured, tested, and its quality is ensured, it can be used in an optimization strategy such as optimal Pareto's optimization.

Then the robustness test proposed here can be applied. The proposed sampling strategy is applied to assemble relevant samples of the evaluated points in the surrogate model throughout the optimization using the meta-heuristic population history. The sampling test is developed so that the collected sample represents the complete population with 99% of statistical confidence. Afterward, the sampling data is assessed as a precise source of information regarding the system and applied to determine the proposed robustness metric. Thereby, it is possible to appraise the surrogate-based optimization's result against a precise source of information. Moreover, computational resources and associated time consumption regarding an optimization based on a rigorous source of information are spared. Additionally, the representation of the analyzed system through the surrogate model optimization is guaranteed with statistical confidence.

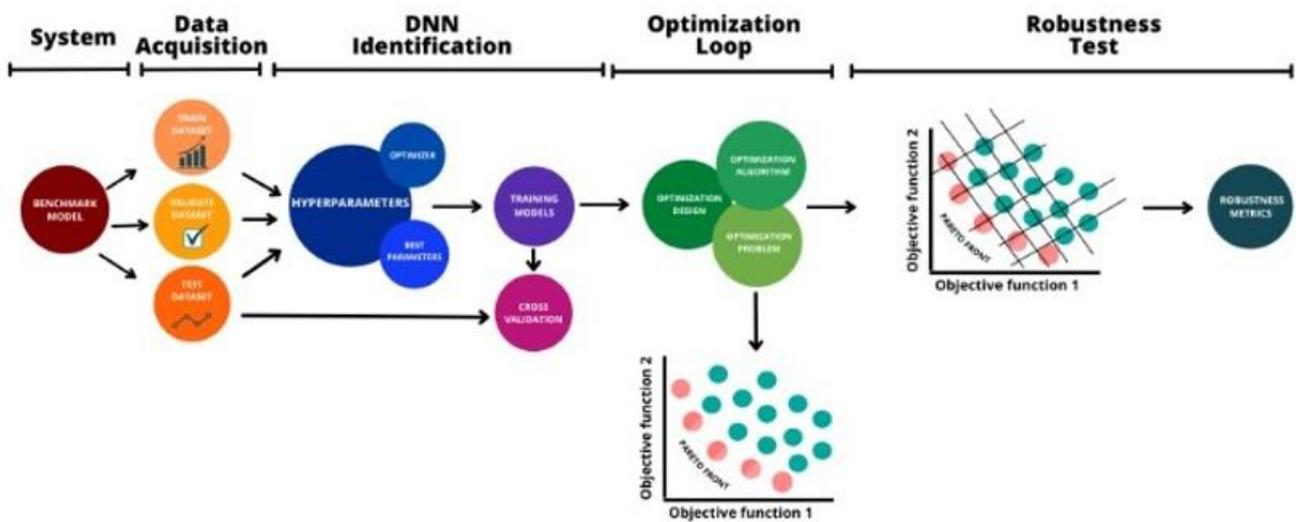

**Figure 1.** Simplified scheme for methodology framework.

*2.1 Design of experiments*

Design of experiments (DoE) is known for bridging the experiments and the model result to improve the optimal results obtained or the quality of a system through statistical methods, such as the Gaussian Process Regression (GPR), the Kiefer-Wolfowitz Equivalence Theorem, and Monte Carlo sampling [9][10]. In this work, a DoE technique is used in two steps. The first step is data acquisition. The DoE is applied in its original context, defining the experimental condition that will allow the SciML models' identification. In the second step of this work, during the proposed robustness analysis, the DoE is applied in identifying reliable samples of the population obtained in the surrogate-based optimization. Thus assuring statistically the representativity of the samples collected, which are used to evaluate the consistency of the SciML model during the optimization procedure.

In this work, the Design of Experiments is proposed based on Latin Hypercube Sampling (LHS). This powerful DoE technique provides a random sampling method that can be useful in sampling with statistical reliability.[11] Olsson et al. (2003) [12] proposed the LHS as a tool to enhance the effectiveness of distinct significant sampling methods applied for structural reliability analysis. Additionally, the mentioned work presented results of 50% of improvement in computer effort savings when comparing the LHS method with the Monte Carlo.

The LHS method is illustrated in Figure 2. This sampling method enables the generation of random samples by dividing the range of each input variable into numerous cells of equal probability. The given input is arranged into layers of probabilistic cells representing its range, which axis is the cumulative probability of the input variables. Afterward, one of the cells is evenly chosen out of the available partitions. Its Cumulative Density Function is locally transposed to obtain random numbers. After that, the cells with similar states or neighboring the chosen one are defined as not available, and the different ones are kept available for the following selection. Then, the indicated procedure keeps the same cycle until it reaches the number of random samples without replacement defined by the user.

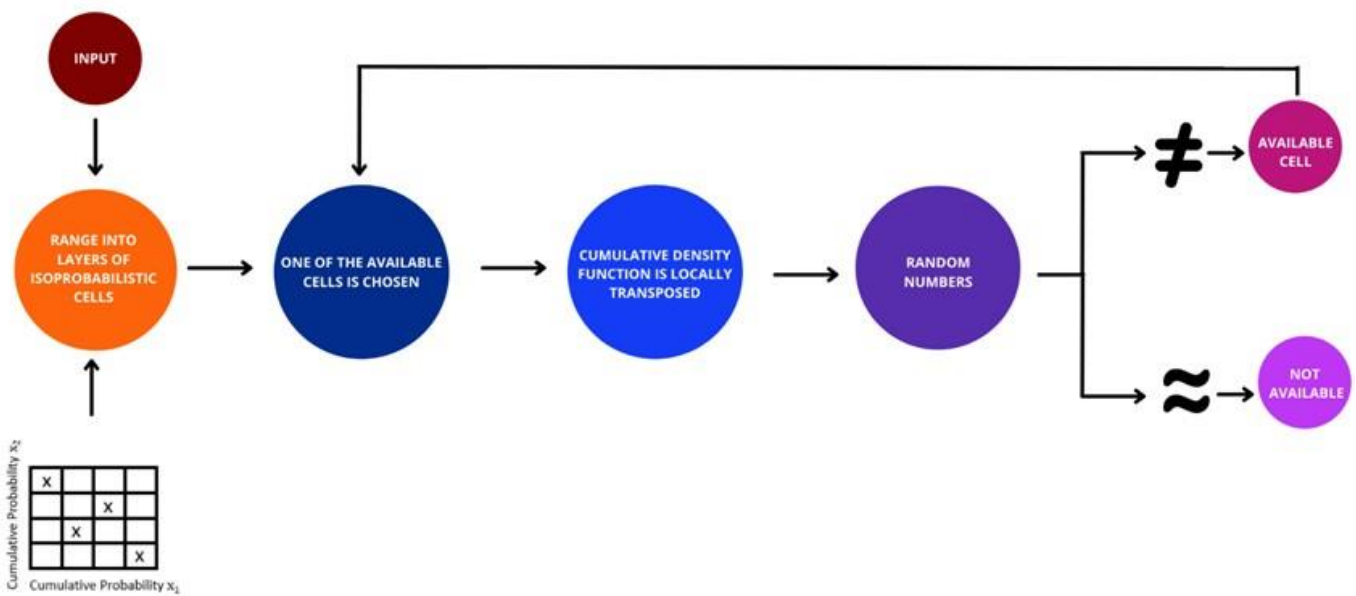

**Figure 2.** Simplified scheme for Latin Hypercube Sampling (LHS).

*2.2 DNN Identification*

The procedure described in section 2.1 is used here to obtain the training, validation, and testing data. After the data acquisition, the SciML model structure must be defined. To accomplish this is necessary to estimate the hyperparameters. These variables dictate the training and performance of a scientific machine learning model. They can be divided into two groups: model and algorithm parameters. The first one contains parameters such as the layer type, number of layers, activation function, and number of neurons. The second group comprises the tuning parameters that can encompass the learning rate, learning rate decay, momentum, and the number of epochs for gradient-based optimizers. The hyperparameters can have different natures, and their estimation is complex. To achieve accurate results, a proper method to search the hyperspace is necessary. The HYPERBAND method is applied in this work, as pointed out in the literature, as a reliable way to identify hyperparameters [13].

The HYPERBAND method was proposed by Li et al. (2018) [13] and is based on appraising the hyperparameter optimization as a pure-exploration problem considering that it applies a present resource allocated to execute random samples of the search space. Furthermore, it is remarkably faster in search speed than the random search and grid methods[13].

Once the model structure is defined by the hyperparameters, the model is again trained to obtain and determine its parameters, weights, and biases. The new training uses all the training, testing, and validation datasets to accomplish the final model structure. Subsequently, the final structure is cross-validated to guarantee the model's quality, precision, and accuracy.

*2.3 Optimization Loop*

In this work, a series of benchmark problems are used to evaluate the proposed methodology, employing it in the context of multiobjective constrained and non-constrained optimizations. Hence, an optimization algorithm is a fundamental part of the system. As this work explores the idea of using a population resulting from a meta-heuristic optimization to provide the robustness analysis of a surrogate-based optimization, a population-based optimization method is necessary.

Eberhart et al. (1995) [14] and Kennedy et al. (1995) [15] developed particle swarm optimization (PSO), a widespread optimization method in the literature. Although it consists of an extraordinarily simple algorithm, PSO was presented as an efficient method for optimizing a broad range of functions and general problems. Rebello et al. (2021) [16] proposed a modified version of the PSO method, the constraints-based sliding particle swarm optimizer (CSPSO). The method presented in the referred work addresses multi-local minima problems. It acquires all minimum points' confidence regions with the entire optimization landscape.

*2.4 Robustness metric*

An optimization procedure will use the SciML models under several conditions, leading it to the edge of its extrapolation capacities. Hence, the SciML model is put under a stress scenario where its robustness is required during the optimization. If not assured, the SciML might lead to the generation of "liars" members in the optimization procedure, consequently leading to deformation in the objective function landscape. In conclusion, in the absence of a robustness test, the results obtained from an optimization procedure based on SciML models should be used with strong parsimony.

Hence, the necessity for an approach where the optimization procedure's robustness can be verified. Thus, allowing broader exploration of such results. The robustness test proposed in this work takes as a base the universal approximation theorem, assuming that a properly identified and robust SciML model will still represent this theorem even in stressful conditions.

**Scientific Machine Learning-based Optimization Robustness Theorem**

Considering the universal approximation theorem, which states that any continuous function such as $F(X, Y)$, $X$ and $Y$ are the set of dependent and independent variables, can be approximated by a Machine Learning model based on an operator, $\delta$, which can be applied to each component of $X$, $\delta(x_i)$, in a way so:

$$\sup_{x \in K} \|f(x) - h(x)\| < \varepsilon \tag{1}$$

With:

$$h(x) = C \cdot \left( \delta(A \cdot x + b) \right)$$

with $\varepsilon$ as a given tolerance, the parameter matrix $A \in \mathbb{R}^{k \times n}$, the bias matrix $b \in \mathbb{R}^{k \times m}$, and the order $k \in \mathbb{N}$. If $n \in \mathbb{N}$, $m \in \mathbb{N}$, then $\delta$ is not a polynomial function, so it is a universal approximator with prediction given by $h(x)$.

The goal is to propose a theorem that demonstrates that after an optimization procedure done using a population-based method, all population members were accurately evaluated by the SciML. Suppose all members of the population respect the theorem presented in equation 1. In that case, the universal approximation is verified, and the SciML model is robust. In this way, a test is sought that is similar to the maximum likelihood, where the probability of $H(X, Y)$ representing $F(X, Y)$ is given by:

$$\phi = P(F(X, Y), H(X, Y), \omega) \tag{2}$$

where $\omega$ is the prediction uncertainty associated with the models' predictions. Hence, if one assumes that:

• The model $F(X, Y)$ is the most rigorus representation of the system;

• The distribution of the deviation between the models can be represented by a gaussian distribution;

• There is no correlation between the successive acquisition of $X$ and $Y$;

• The uncertainties are negligible.

Then, the likelihood of $H(X, Y)$ representing $F(X, Y)$ for a specific prediction made by the SciML model can be given by:

$$\phi_i = \frac{1}{\sqrt{2\pi\sigma^2}} e^{-\frac{\left(f(x_i, y_i) - h(x_i, y_i)\right)^2}{2\sigma^2}} \tag{3}$$

where $\sigma$ is the variance. From equation 3, it is possible to extend the likelihood function toward all predictions, which can be written by:

$$\phi = \prod_{i=1}^{Np} \frac{1}{\sqrt{2\pi\sigma^2}} e^{-\frac{\left(f(x_i, y_i) - h(x_i, y_i)\right)^2}{2\sigma^2}} \tag{4}$$

where $Np$ is the total number of predictions done by the SciML model during the optimization procedure. However, the above function is difficult to be differentiated. In order to overcome this problem, we can apply logarithm properties in equation 4, obtaining:

$$\phi = Np \ln \frac{1}{\sqrt{2\pi\sigma^2}} - \frac{1}{2\sigma^2} \sum_{i=1}^{Np} \left(f(x_i, y_i) - h(x_i, y_i,)\right)^2 \tag{5}$$

Finally, differentiating equation 5 as a function of $\delta$, one can obtain the equation below:

$$\frac{d\phi}{d\delta} = \frac{1}{2\sigma^2} \sum_{i=1}^{Np} 2\big(f(x_i, y_i) - h(x_i, y_i,)\big) \tag{6}$$

Rewriting equation 6, one obtains:

$$\frac{d\phi}{d\delta} = \frac{2}{2\sigma^2} \sum_{i=1}^{Np} \big(f(x_i, y_i) - h(x_i, y_i,)\big) \tag{7}$$

In such a case, the maximum likelihood will lead to a point where the above derivative equals zero. Consequently, the likelihood test respects the universal approximation theorem, equation 1. Hence, the final form of the robustness test can be written as:

$$\sum_{i=1}^{Np} \big(f(x_i, y_i) - h(x_i, y_i,)\big) = 0 \tag{8}$$

Then, a robustness test is deduced that guarantees the universal approximation theorem, equation 1. As it is possible to see, we reached an equation that resembles the theorem presented in equation 1. The robustness metric ($Rb$) can be computed as:

$$Rb = \sum_{i=1}^{Np} \big(f(x_i, y_i) - h(x_i, y_i,)\big) \tag{9}$$

The test deduced considers the existence of a rigorous source of information regarding the system. From equations 1 and 8, the test can be written as:

$$Rb \leq \varepsilon \tag{10}$$

Then, as close as the robustness metric ($Rb$) is to $\varepsilon$, as robust is the optimization result. This can be either a rigorous model or a data set observed in the real system. The deduced test presupposes the evaluation of all members of the optimization papulation using a rigorous source of information. This would lead to a critical limitation of the proposed methodology. However, as it is a population-based approach, one can apply a statistical test to evaluate the minimum size of the population that can be drawn to have a high level of statistical confidence that the sample represents the entire population. If the total size of the population is known, this can be computed as:

$$Np = \frac{(\text{Zscore})^2 \, \sigma \, (1-\sigma)}{(\text{E})^2} \qquad (11)$$

Zscore, the critical value corresponding to the desired confidence level, $E$ is the maximum allowed difference between the population average and sampled population average.

The final part of the proposed methodology assesses the robustness of the optimal result achieved. The procedure established here to certify the result is well-grounded on a meta-heuristic method. The optimization points are randomly sampled in their population through the Latin Hypercube Sampling (LHS) previously explained. From the random and non-correlated sampled population, their corresponding states are evaluated in a rigorous source of information, and the robustness metric is computed.

## 3. Discussion

### *3.1 Benchmark 1 – Binh and Korn function*

The first benchmark function tested is the Binh and Korn function, a convex, nonlinear, constrained, two-parameter function. The Binh and Korn function benchmark consists of two objective functions subjected to two constraints described as follows:

$$\text{Minimize } F = (F_1(x_1, x_2), \ F_2(x_1, x_2)) \tag{12}$$

$$F_1(x_1, x_2) = 4x_1{}^2 + 4x_2{}^2 \tag{13}$$

$$F_2(x_1, x_2) = (x_1 - 5)^2 + (x_2 - 5)^2 \tag{14}$$

where

$$G_1(x_1, x_2) = (x_1 - 5)^2 + x_2{}^2 - 25 \leq 0 \tag{15}$$

$$G_2(x_1, x_2) = -(x_1 - 8)^2 - (x_2 + 3)^2 + 7.7 \leq 0 \tag{16}$$

$$0 \leq x_1 \leq 5 \tag{17}$$

$$0 \leq x_2 \leq 3 \tag{18}$$

The LHS is then used to generate random synthetic data to identify the surrogate model. As exemplified in Figure 3, there is no correlation between the variables in any of the objective functions or constraints. For that reason, the LHS application can be considered consistent for the data generation and sampling steps.

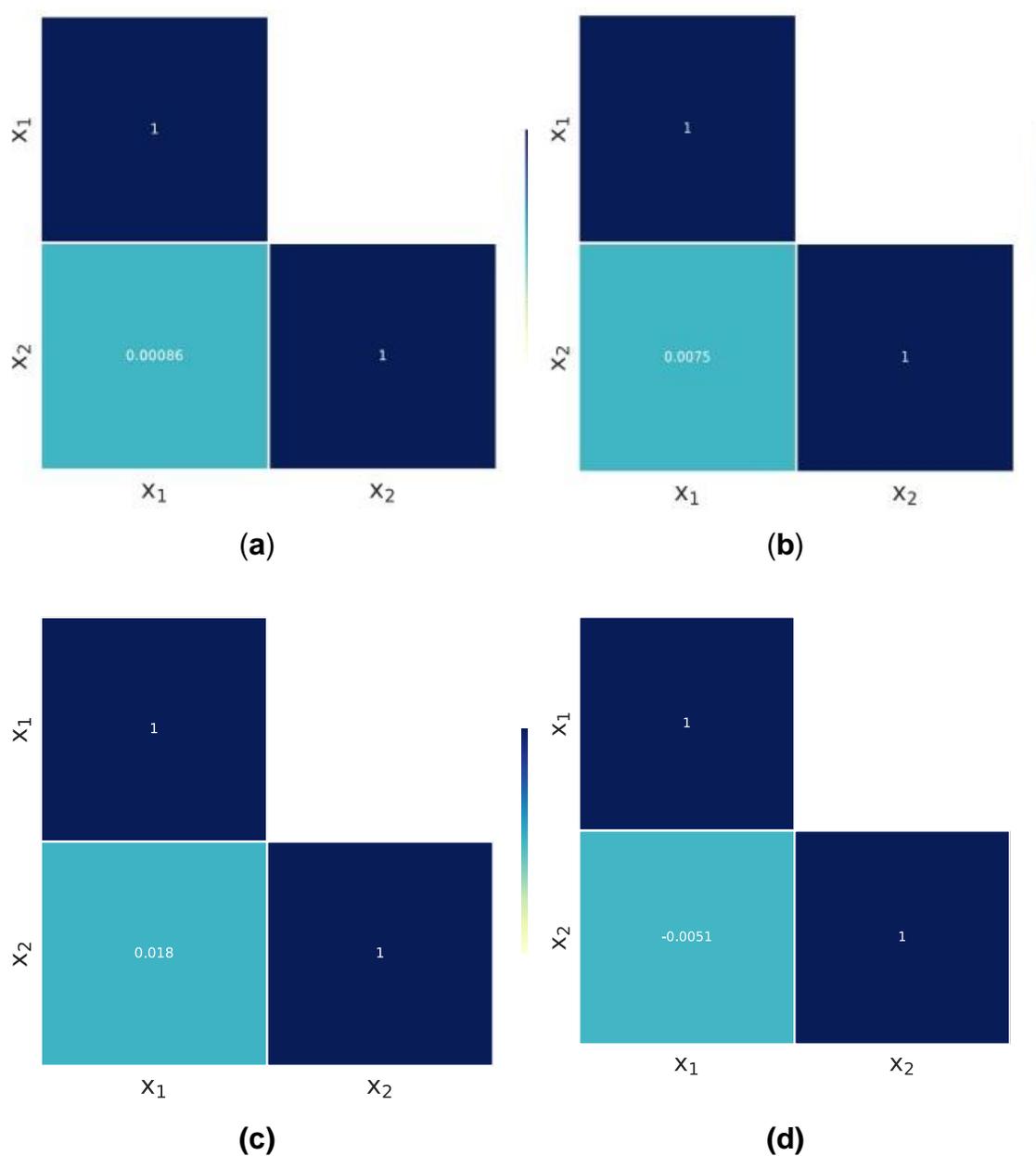

**Figure 3.** Correlation Matrix for Binh and Korn function: (**a**)Objective function 1; (**b**)Objective function 2; (**c**) Constrain 1; (**d**)Constrain 2.

Then, following the proposed methodology, the SciML models are identified. The previously described methodology is here applied to identify a surrogate model that can represent properly $F_1$ and $F_2$ in a MISO strategy (Multi-Input Single Output). In this case, a model for each variable that represents these functions is generated, validated, and tested. The best model to represent the testing benchmark 1 has the structure depicted in Table 1.

**Table 1.** Best model structure based on hyperparameters for Binh and Korn function.

| Objective | Initial learning rate | Batch size | Number of neurons Intermediate fully connected layer | Number of parameters Intermediate fully connected layer | Total number of parameters | Number of trainable parameters |
|---|---|---|---|---|---|---|
| Objective function 1 | 0.001 | 32 | {Layer 1: 60, Layer 2: 60, Layer 3: 60, Layer 4: 1} | {Layer 1: 180, Layer 2: 3660, Layer 3: 3660, Layer 4: 61} | 7561 | 7561 |
| Objective function 2 | 0.001 | 32 | {Layer 1: 60, Layer 2: 60, Layer 3: 60, Layer 4: 1} | {Layer 1: 180, Layer 2: 3660, Layer 3: 3660, Layer 4: 61} | 7561 | 7561 |
| Constrain 1 | 0.001 | 32 | {Layer 1: 60, Layer 2: 60, Layer 3: 60, Layer 4: 1} | {Layer 1: 180, Layer 2: 3660, Layer 3: 3660, Layer 4: 61} | 7561 | 7561 |
| Constrain 2 | 0.001 | 32 | {Layer 1: 60, Layer 2: 60, Layer 3: 60, Layer 4: 1} | {Layer 1: 180, Layer 2: 3660, Layer 3: 3660, Layer 4: 61} | 7561 | 7561 |

Figure 4 illustrates the comparison of the surrogate model's prevision and the benchmark's response; it shows a good agreement and precision of the surrogate model with the synthetic data. This conclusion can be reinforced by the parity chart that exhibits the prediction distribution related to the test data, Figure 5. In this case, the data is entirely centered and randomly dispersed in the sloping curve, indicating that the model can precisely represent the benchmark.

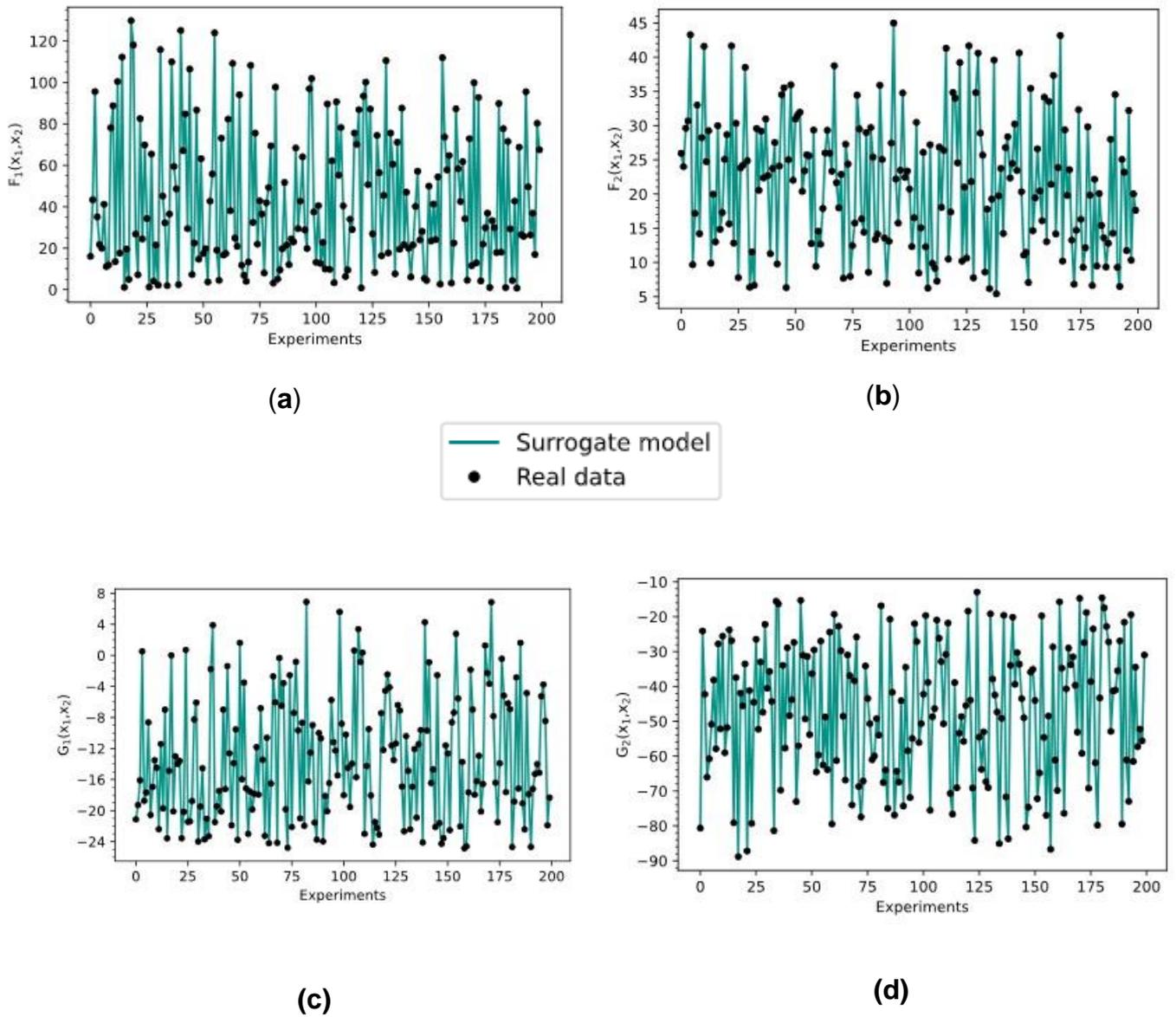

**Figure 4.** Comparison graph for Binh and Korn function real data and surrogate model representation: (**a**)Objective function 1; (**b**)Objective function 2; (**c**) Constrain 1; (**d**)Constrain 2.

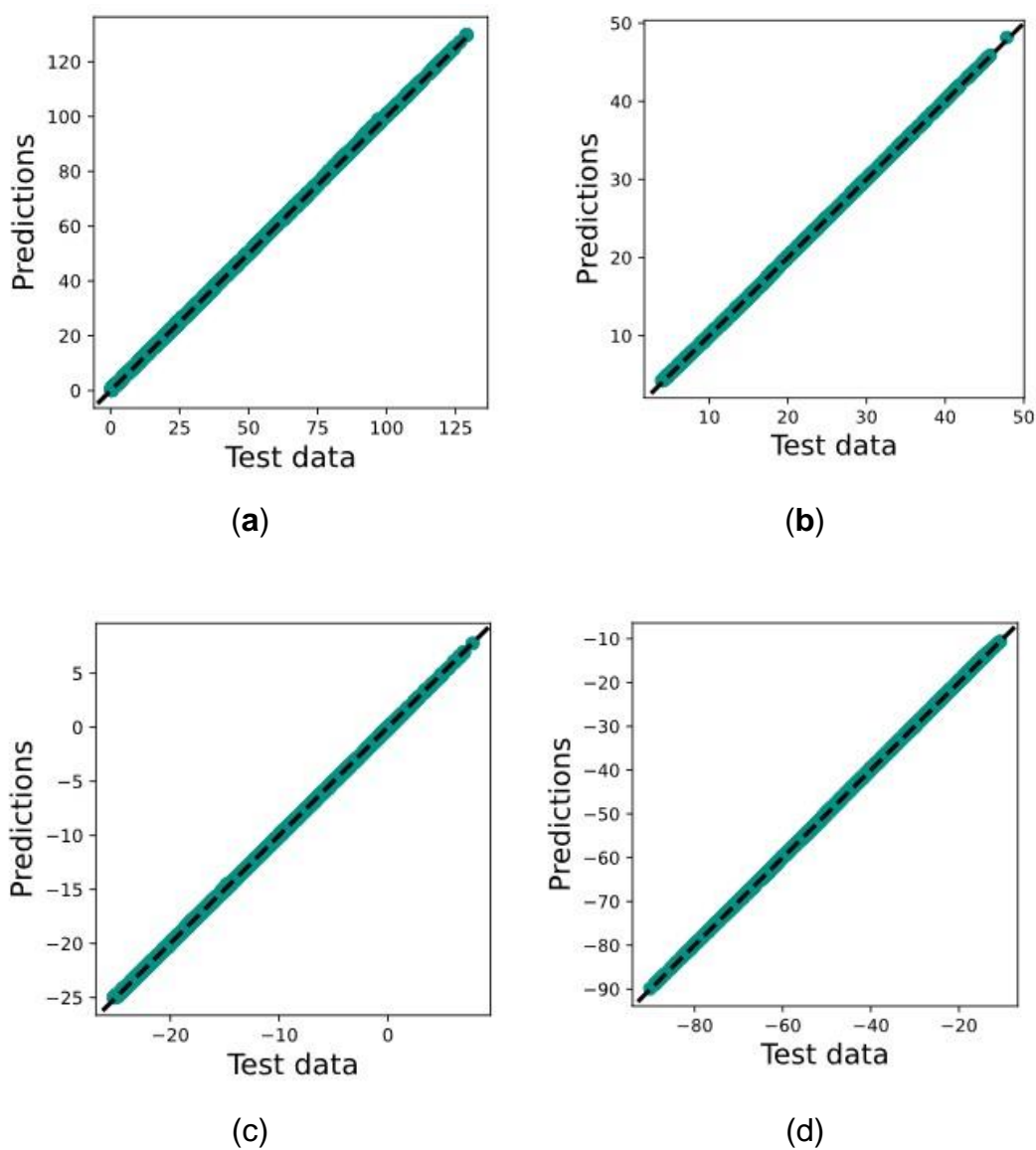

**Figure 5.** Parity chart for Binh and Korn function: (**a**)Objective function 1; (**b**)Objective function 2; (**c**) Constrain 1; (**d**) Constrain 2.

Since the models are already tested, validated, and their efficiency and precision are evaluated, $F_1$ and $F_2$ can be represented with high precision, and the CSPSO optimizer can be applied to solve the previously mentioned optimization problems. Mean Squared Error (M.S.E.) and Mean Absolute Error (M.A.E.) values are obtained for each objective function and their constraints, Table 2.

**Table 2.** Mean Squared Error (MSE) and Mean Absolute Error (MAE) values for Binh and Korn function.

| Objective | MSE | MAE |
|---|---|---|
| Objective function 1 | 0.0000076872 | 0.0021 |
| Objective function 2 | 0.000003728 | 0.0015 |
| Constrain 1 | 0.0000045345 | 0.0017 |
| Constrain 2 | 0.0000052308 | 0.002 |

Figure 6 contains the resulting feasible regions and Pareto fronts obtained from the CSPSO optimization. The Binh and Korn function was optimized through rigorous models, aiming to test and validate the obtained result. This optimization result is also presented in Figure 6. It is possible to visualize an overlap of the Pareto fronts obtained for both optimizations. These results attest to the consistency of the methodology applied and infer that surrogate-based models can be confidently used to solve optimization problems as they provide precise results. However, repeating the rigorous model's optimization is exhaustive and not feasible in real applications. Therefore, the robustness test is proposed here to provide a viable and reliable way to test the consistency of such results.

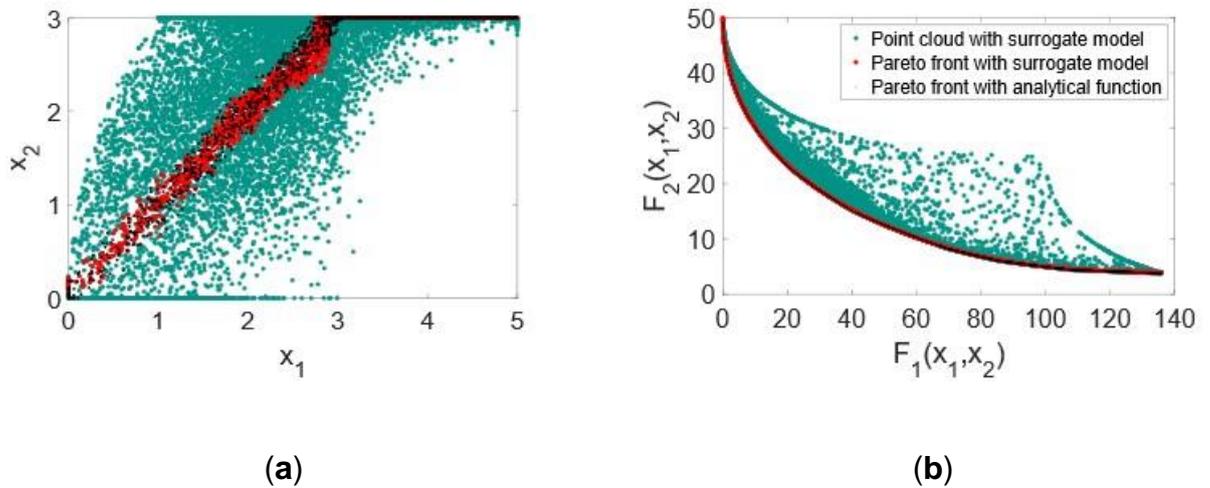

(**a**)                                             (**b**)

**Figure 6.** Feasible regions and Pareto fronts of rigorous model and CSPSO model: (**a**) Decision variables; (**b**) Objective functions.

In closing, a robustness test is applied to the results obtained for the surrogate-based optimization of the benchmark. For this instance, the LHS is again used to withdraw random samples without replacement from the feasible region of the optimization population obtained through the CSPSO and assemble significant representatives of the region's topography from the objective functions. Following the proposed methodology, the optimization's robustness is guaranteed, only requiring a small percentage of the resulting particles from the CSPSO for the evaluation. Equation 11 is applied to compute the minimum number of samples needed to obtain the confidence of 99% that it represents the total population. Then, the $Np$ was computed, which was equal to 381 samples in this case. These were withdrawn to enable an expeditious evaluation of the result's robustness without needing a new optimization using the rigorous model.

The values of $Rb$ and M.A.E. for this test are presented in table 3. Considering a $\varepsilon = 0.05$, the Rb values under the $\varepsilon$ threshold affirm the precision associated with the surrogate-based multiobjective optimization. It is important to note that this test reduces from 10000 evaluations of the objective function necessary to assure that the optimization is robust to only 381.

**Table 3.** Robustness metric and Mean Absolute Error (M.A.E.) values for the robustness test.

| Objective | Rb | MAE |
|---|---|---|
| Objective function 1 | 0.0214 | 0.0008 |
| Objective function 2 | 0.0648 | 0.0063 |

### 3.2 Benchmark 2 – Zitzler–Deb–Thiele's function N. 3

The Zitzler–Deb–Thiele's function N. 3 is a box-constrained, continuous multiobjective problem with up to thirty dimensions. In this case, this benchmark consists of two objective functions described by three variables as shown below:

$$\text{Minimize } F = (F_1(x_1) \, , \; F_2(x_1, x_2, x_3)) \tag{19}$$

$$F_1(x_1) = \; x_1 \tag{20}$$

$$F_2(x_1, x_2, x_3) = 1 + \frac{9}{29}\sum_{i=2}^{3} x_i \; \left(1 - \sqrt{\frac{x_1}{\frac{9}{29}\sum_{i=2}^{3} x_i}} - \left(\frac{x_1}{\frac{9}{29}\sum_{i=2}^{3} x_i}\right)\sin(10\pi x_1)\right) \tag{21}$$

where

$$0 \leq x_1, x_2, x_3 \leq \; 1 \tag{22}$$

Figure 7 attests that there is no correlation between the variables in any of the objective functions. The LHS application for the Zitzler–Deb–Thiele's function N. 3 can be reliable for data generation.

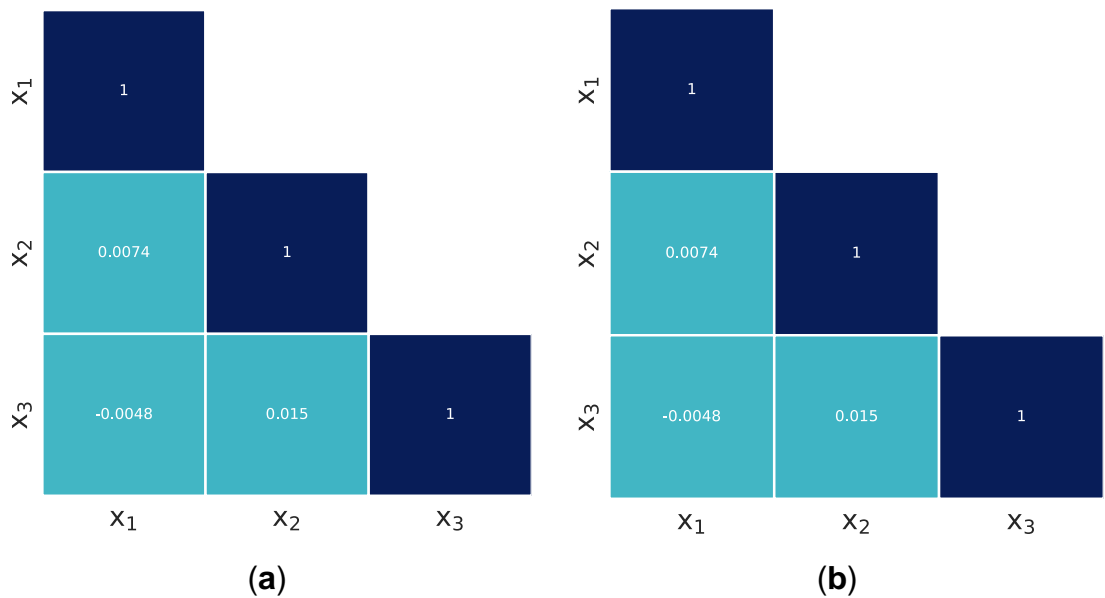

(**a**)                                        (**b**)

**Figure 7.** Correlation Matrix for Zitzler–Deb–Thiele's function N. 3: (**a**) Objective function 1; (**b**) Objective function 2.

The assessment of the second benchmark follows the previously proposed methodology to identify the surrogate model that is more consistent in the representation of $F_1$ and $F_2$, applying the MISO strategy for this instance. Applying the previously presented methodology resulted in the SciML models' structure shown in Table 4 as the best model representing the testing benchmark 2.

**Table 4.** Best model structure based on hyperparameters for Zitzler–Deb–Thiele's function N. 3.

| Objective | Initial learning rate | Batch size | Number of neurons Intermediate fully connected layer | Number of parameters Intermediate fully connected layer | Total number of parameters |
|---|---|---|---|---|---|
| Objective function 1 | 0.001 | 32 | {Layer 1: 100, Layer 2: 1} | {Layer 1: 400, Layer 2: 101} | 501 |
| Objective function 2 | 0.001 | 32 | {Layer 1: 80, Layer 2: 80, Layer 3: 80, Layer 4: 80, Layer 5: 80, Layer 6: 1} | {Layer 1: 320, Layer 2: 6480, Layer 3: 6480, Layer 4: 6480, Layer 5: 6480, Layer 6: 81} | 26321 |

Figure 8 presents the surrogate model's good agreement and precision with the synthetic test data. Figure 9 demonstrates that the model can efficiently represent the benchmark.

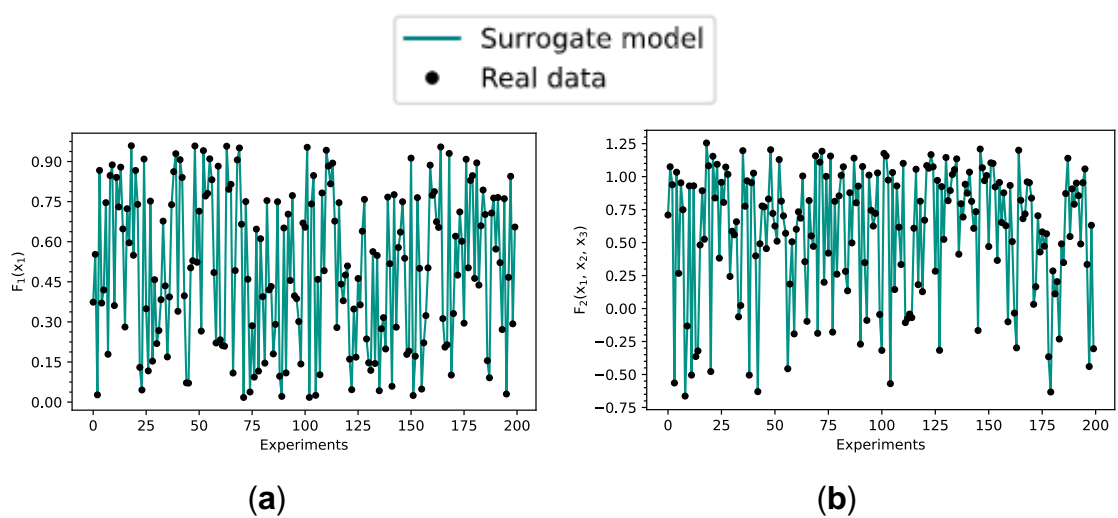

(**a**)          (**b**)

**Figure 8.** Comparison graph for Zitzler–Deb–Thiele's function N. 3 real data and surrogate model representation: (**a**)Objective function 1; (**b**)Objective function 2.

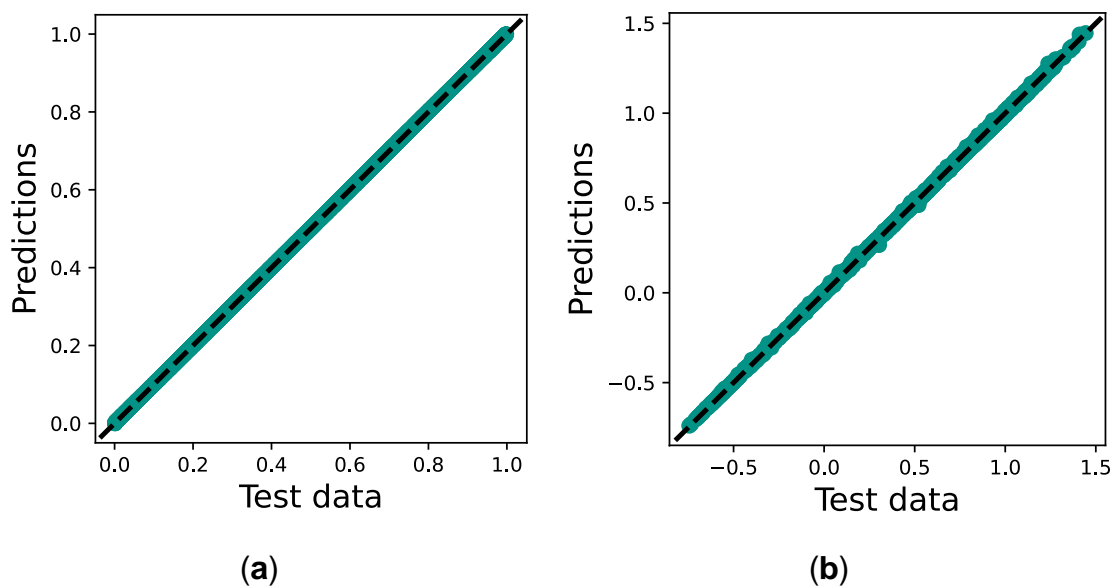

(**a**)                              (**b**)

**Figure 9.** Parity chart for Zitzler–Deb–Thiele's function N. 3: (**a**)Objective function 1; (**b**)Objective function 2.

The results for M.S.E. and M.A.E. obtained after validation, testing and evaluation of the models are presented in Table 5.

**Table 5.** Mean Squared Error (M.S.E.) and Mean Absolute Error (M.A.E.) values for Zitzler–Deb–Thiele's function N. 3.

| Objective | MSE | MAE |
|---|---|---|
| Objective function 1 | 0.0000011359 | 0.001 |
| Objective function 2 | 0.0000059264 | 0.002 |

The resulting feasible regions and Pareto fronts from the implementation of the CSPSO model and the rigorous model's optimization are presented in Figure 10. The graphs show an overlap of the Pareto fronts obtained for both optimizations, which certify the consistency of the methodology applied and make it possible to conclude that surrogate-based models applied to solve optimization problems can provide precise results.

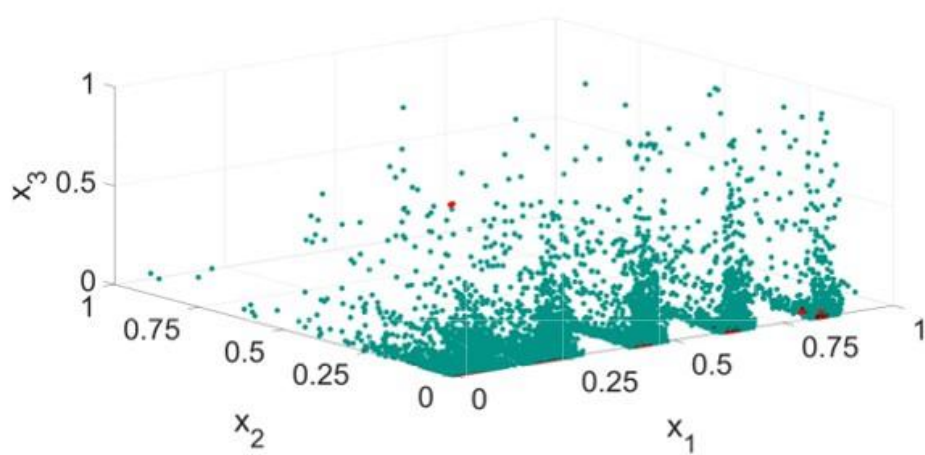

**(a)**

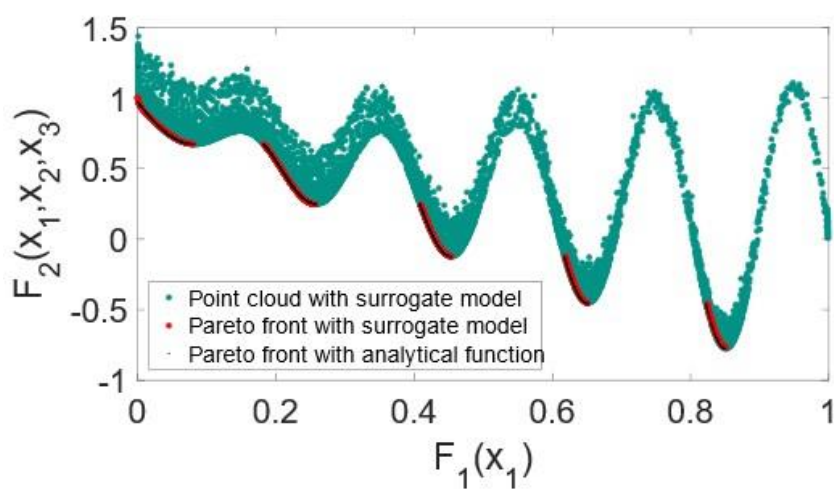

(**b**)

**Figure 10.** Feasible regions and Pareto fronts of rigorous model and CSPSO model: (**a**) Decision variables; (**b**) Objective functions.

The robustness test is performed analogously to the previous benchmark analyzed to evaluate the results achieved for the surrogate-based model application. In resemblance, through Equation 10 it was possible to obtain the required number of samples to ensure the 99% confidence in the representation of the total population, and 382 samples were acquired.

Following the same path, Table 6 presents the values for Rb and M.A.E. for the present robustness test. The Rb values far below the 0.05 threshold confirm the high precision of the surrogate model's representation of the system.

**Table 6.** Robustness metric and Mean Absolute Error (M.A.E.) values for the robustness test.

| Objective | Rb | MAE |
|-----------|-----|-----|
| Objective function 1 | 0.00485 | 0.0000564 |
| Objective function 2 | 0.00655 | 0.000113 |

### 3.3 Benchmark 3 – DTLZ2

The DTLZ2 function is continuous in the search space and unimodal and consists of three objective functions and one box constraint. The final representation is a 3D space. The functions that are comprehended in the benchmark are written below:

$$\text{Minimize } F = (F_1(x_1, x_2, x_3), F_2(x_1, x_2), \ F_3(x_1)) \tag{23}$$

$$F_1(x_1, x_2, x_3) = (1+G(x_M))\cos(\tfrac{x_1\pi}{2})\cos(\tfrac{x_2\pi}{2})\sin(\tfrac{x_3\pi}{2}) \tag{24}$$

$$F_2(x_1, x_2) = (1+g(x_M))\cos(\tfrac{x_1\pi}{2})\sin(\tfrac{x_2\pi}{2}) \tag{25}$$

$$F_3(x_1) = (1+G(x_M))\sin(\tfrac{x_1\pi}{2}) \tag{26}$$

where

$$G(x_1, x_2, x_3) = \sum_{i=1}^{3}(x_i - 0.5)^2 \tag{27}$$

$$0 \leq x_1, x_2, x_3 \leq 1 \tag{28}$$

Correspondingly to the procedure implemented in the previous section, the LHS method was once again applied to generate random experimental data to evaluate the surrogate model. Figure 11 attests that there is no correlation between the variables in any of the objective functions and that the LHS application for the DTLZ2 is reliable.

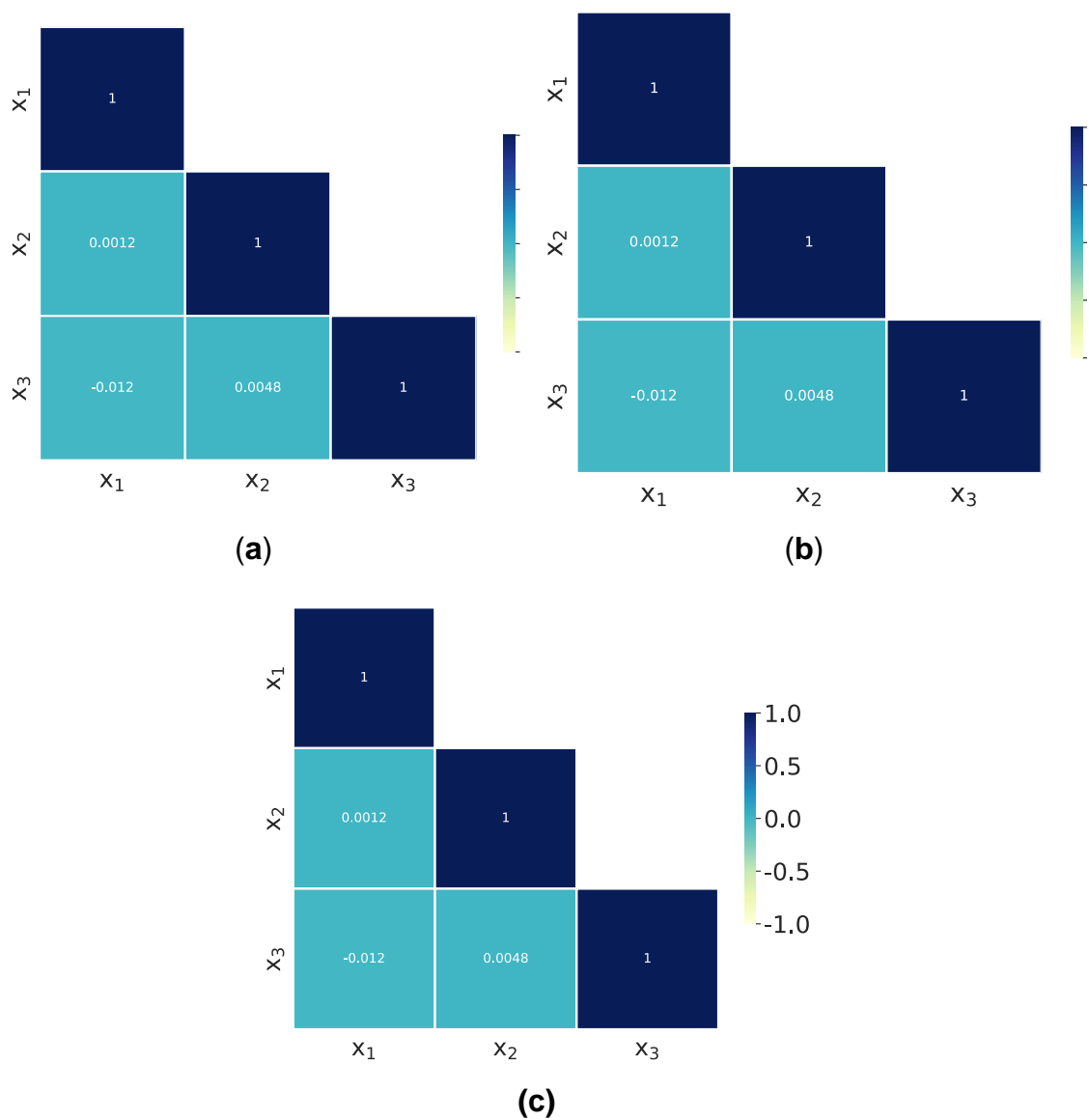

**Figure 11.** Correlation Matrix for DTLZ2 function: (**a**) Objective function 1; (**b**) Objective function 2; (**c**) Objective function 3.

Similarly, the proposed methodology is applied to the framework followed in the other two benchmarks. The SciML model structure is presented in Table 7 as the best model to represent benchmark 3. It is sought to identify the surrogate model that is more consistent in the representation of $F_1$, $F_2$ and $F_3$. In this context, the MISO strategy is applied.

**Table 7.** Best model structure based on hyperparameters for DTLZ2 function.

| Objective | Initial learning rate | Batch size | Number of neurons Intermediate fully connected layer | Number of parameters Intermediate fully connected layer | Total number of parameters |
|---|---|---|---|---|---|
| Objective function 1 | 0.001 | 32 | {Layer 1: 180, Layer 2: 180, Layer 3: 180, Layer 4: 180, Layer 5: 180, Layer 6: 1} | {Layer 1: 720, Layer 2: 32580, Layer 3: 32580, Layer 4: 32580, Layer 5: 32580, Layer 6: 181} | 131221 |
| Objective function 2 | 0.001 | 32 | {Layer 1: 220, Layer 2: 220, Layer 3: 220, Layer 4: 220, Layer 5: 220, Layer 6: 1} | {Layer 1: 880, Layer 2: 48620, Layer 3: 48620, Layer 4: 48620, Layer 5: 48620, Layer 6: 221} | 195581 |
| Objective function 3 | 0.001 | 32 | {Layer 1: 220, Layer 2: 220, Layer 3: 1} | {Layer 1: 880, Layer 2: 48620, Layer 3: 221} | 49721 |

The third benchmark is tested similarly to the previous ones. Figure 12 illustrates the surrogate model's good agreement and precision with the synthetic data. This result is strengthened by Figure 13 demonstrating that the model can efficiently represent the benchmark.

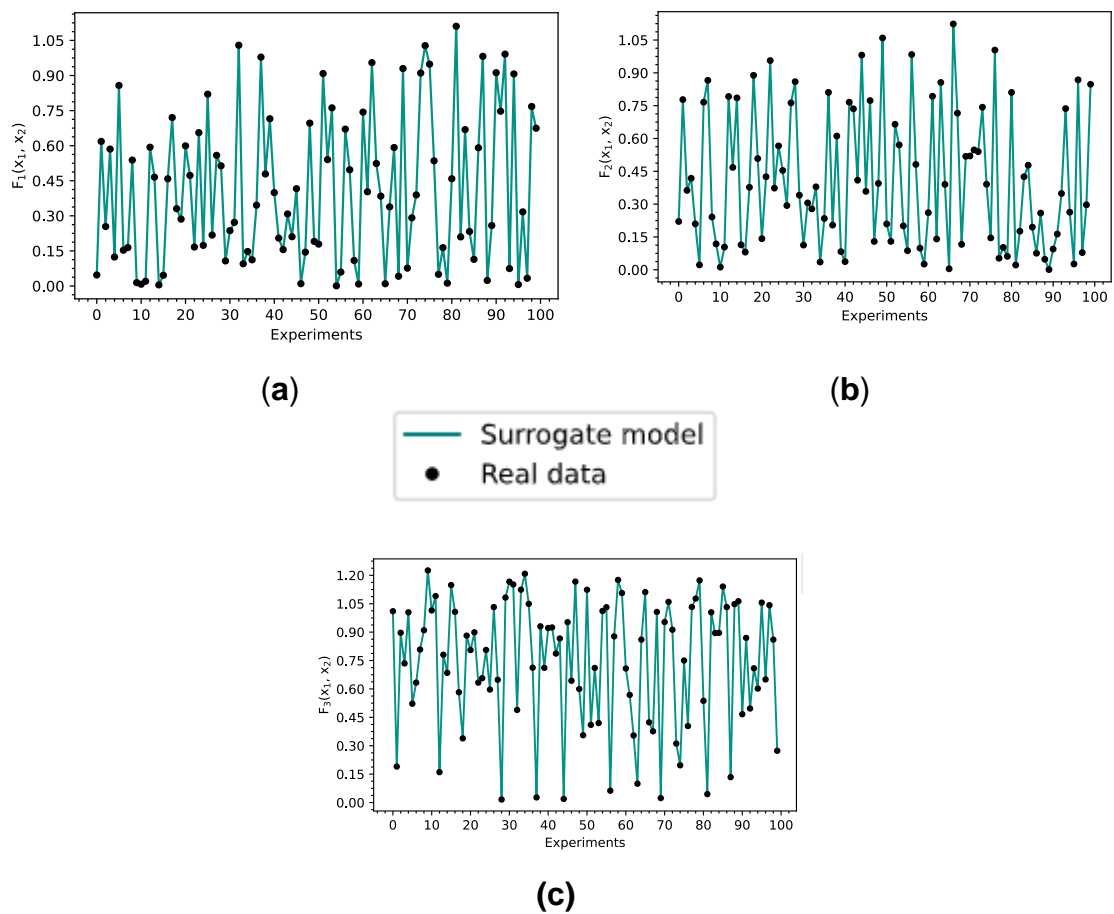

**Figure 12.** Comparison graph for DTLZ2 real data and surrogate model representation: (**a**) Objective function 1; (**b**) Objective function 2; (**c**) Objective function 3.

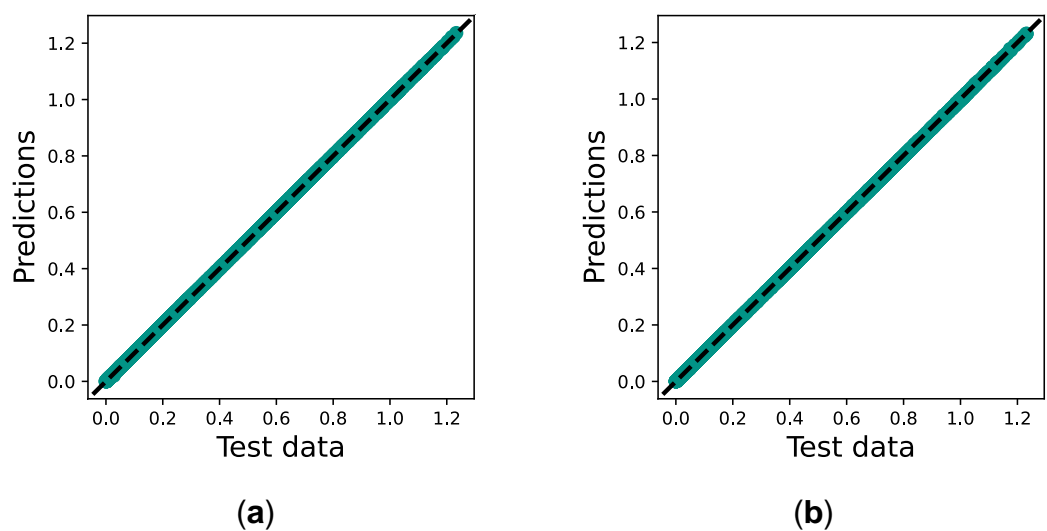

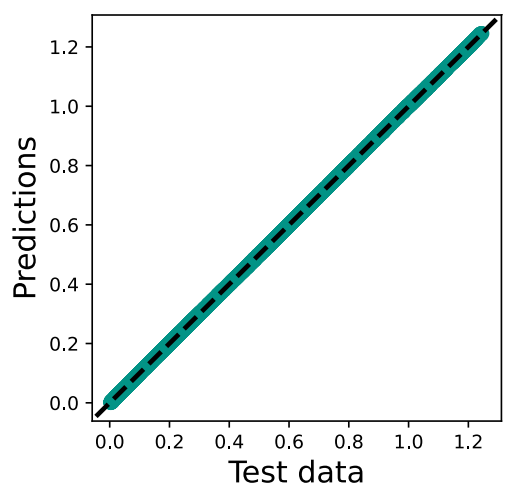

**(c)**

**Figure 13.** Parity chart for DTLZ2 function: (**a**) Objective function 1; (**b**) Objective function 2; (**c**) Objective function 3.

The models are then validated, tested, and evaluated. After that, results for M.S.E. and M.A.E. are obtained and presented in Table 8.

**Table 8.** Mean Squared Error (MSE) and Mean Absolute Error (MAE) values for DTLZ2 function.

| Objective | MSE | MAE |
|---|---|---|
| Objective function 1 | 0.00000072915 | 0.00068414 |
| Objective function 2 | 0.00000061084 | 0.00059728 |
| Objective function 3 | 0.00000085066 | 0.00073728 |

The feasible regions and Pareto fronts resulting from the implementation of the CSPSO model and the rigorous model's optimization are presented in Figure 14. The graphs show an overlap of the Pareto fronts obtained for both optimizations, demonstrating that the methodology applied is consistent and the surro-gate-based model applied to solve the proposed optimization problem can provide precise results. However, this demonstration comes with a cost. It was necessary to reevaluate the objective function another 10000 times.

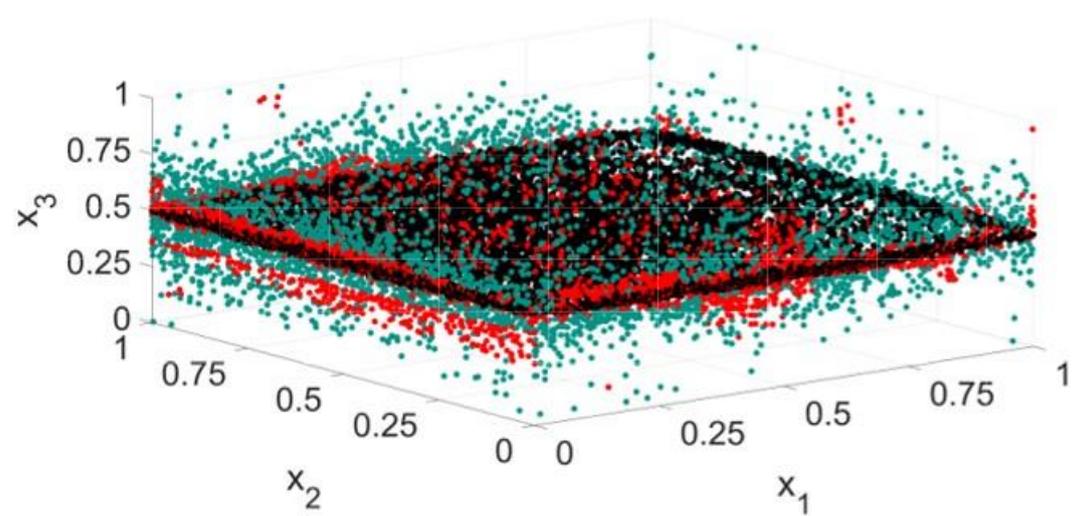

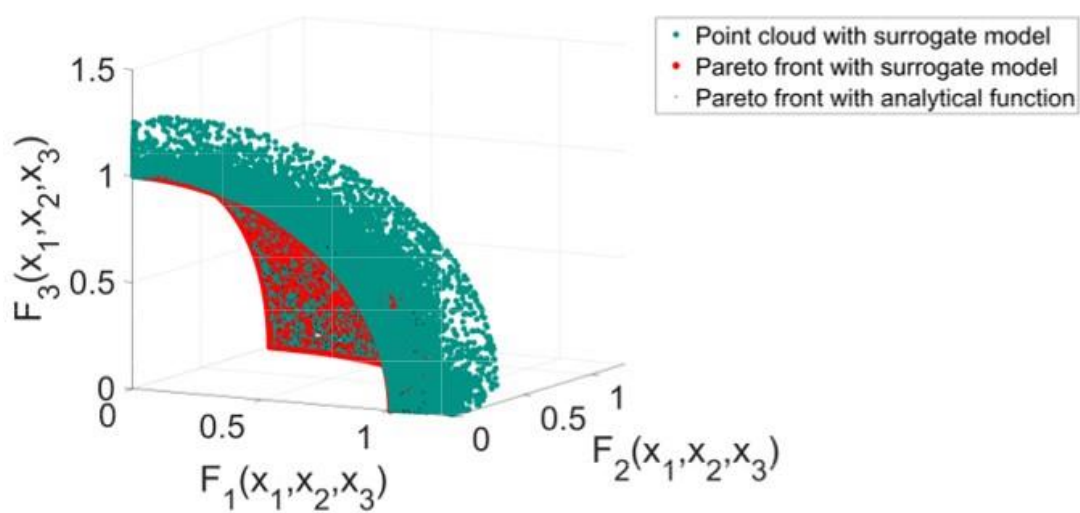

(**b**)

**Figure 14.** Feasible regions and Pareto fronts of rigorous model and CSPSO model: (**a**) Decision variables; (**b**) Objective functions.

The surrogate model's result is then assessed in terms of robustness, as proposed here. To accomplish this test, the proposed methodology is applied. Equation 11 is then used to obtain the required number of samples to achieve 99% confidence of representativity. The value of 383 samples was obtained. The resultant values for Rb and M.A.E. are presented in Table 9. Their reduced values reiterate the model's precision, fitting, and robustness validation. Again, the Rb metric is far below the $\varepsilon$, of 0.05.

**Table 9.** Robustness metric and Mean Absolute Error (M.A.E.) values for the robustness test.

| Objective | Rb | MAE |
|---|---|---|
| Objective function 1 | 0.000911 | 0.00000132 |
| Objective function 2 | 0.000717 | 0.00000084 |
| Objective function 3 | 0.001305 | 0.00000283 |

## 5. Conclusions

This work presents a novel strategy to guarantee the robustness of Scientific Machine Learning-based optimization. The surrogate models offer remarkable advantages for optimization, such as the reduced computation resources required and the potential of application for real-time optimization. However, their robustness is a concern when applying this approach in optimization frameworks.

This work deduced a robustness test to guarantee that the optimization results respect the universal approximation theorem. Then the test is used within the proposed framework. The proposed methodology is evaluated through several benchmark multiobjective optimization tests. The robustness test of the surrogate model optimization is proposed. Given reliability by a 99% confidence regarding the representation of the total population, it is possible to verify that the surrogate-based model can be used with robustness within the optimization loop. The results are robust and well-fitted for representing real scenarios. Therefore, following the framework methodology described in this work it is possible to obtain a robust optimization with a high precision level.

## Acknowledgments

This work was financially supported by LA/P/0045/2020 (ALiCE), UIDB/50020/2020 and UIDP/50020/2020 (LSRE-LCM), funded by national funds through FCT/MCTES (PIDDAC), by and F.C.T.—Fundação para a Ciência e Tecnologia under C.E.E.C. Institucional program. Capes for its financial support, financial code 001